\documentclass[12pt]{amsart}
\usepackage{thmtools}
\usepackage{amsfonts,amsmath,amssymb}
\usepackage{hyperref}
\usepackage{paralist}
\usepackage{stmaryrd}
\usepackage[linesnumbered,ruled,vlined,norelsize]{algorithm2e}
\usepackage[utf8]{inputenc}

\declaretheoremstyle[bodyfont=\normalfont]{noncursive}
\declaretheorem{theorem}
\declaretheorem[numberwithin=section]{lemma}

\declaretheorem[numberlike=lemma]{proposition}
\declaretheorem[numberlike=lemma]{corollary}
\declaretheorem[style=noncursive,numberlike=lemma]{definition}

\declaretheorem[style=noncursive,numberlike=lemma]{remark}

\usepackage{amsmath,amssymb,amsthm,enumerate}

\begin{document}

\numberwithin{equation}{section}

\def\1#1{\overline{#1}}
\def\2#1{\widetilde{#1}}
\def\3#1{\widehat{#1}}
\def\4#1{\mathbb{#1}}
\def\5#1{\frak{#1}}
\def\6#1{{\mathcal{#1}}}

\def\C{{\4C}}
\def\R{{\4R}}
\def\N{{\4N}}
\def\Z{{\4Z}}

\def \im{\text{\rm Im }}
\def \re{\text{\rm Re }}
\def \Char{\text{\rm Char }}
\def \supp{\text{\rm supp }}
\def \codim{\text{\rm codim }}
\def \Ht{\text{\rm ht }}
\def \Dt{\text{\rm dt }}
\def \hO{\widehat{\mathcal O}}
\def \cl{\text{\rm cl }}
\def \bR{\mathbb R}
\def \bC{\mathbb C}
\def \bP{\mathbb P}
\def \C{\mathbb C}
\def \bL{\mathbb L}
\def \bZ{\mathbb Z}
\def \bN{\mathbb N}
\def \scrF{\mathcal F}
\def \scrK{\mathcal K}
\def \scrM{\mathcal M}
\def \cR{\mathcal R}
\def \scrJ{\mathcal J}
\def \scrA{\mathcal A}
\def \scrO{\mathcal O}
\def \scrV{\mathcal V}
\def \scrL{\mathcal L}
\def \scrE{\mathcal E}

\oddsidemargin=0.1in \evensidemargin=0.1in \textwidth=6.4in
\headheight=.2in \headsep=0.1in \textheight=8.4in
\newcommand{\rl}{{\mathbb{R}}}
\newcommand{\cx}{{\mathbb{C}}}
\newcommand{\id}{{\mathbb{I}}}
\newcommand{\m}{{\mathcal{M}}}
\newcommand{\dbar}{\overline{\partial}}
\newcommand{\Db}[1]{\frac{\partial{#1}}{\partial\overline{z}}}
\newcommand{\abs}[1]{\left|{#1}\right|}
\newcommand{\e}{\varepsilon}
\newcommand{\tmop}[1]{\ensuremath{\operatorname{#1}}}
\renewcommand{\Re}{\tmop{Re}}
\renewcommand{\Im}{\tmop{Im}}
\newcommand{\dist}{{\mathrm{dist}}}
\newcommand {\OO}{{\mathcal O}}
\newcommand {\Sig}{{\Sigma}}
\renewcommand {\a}{\alpha}
\newcommand {\MC}{M^{\mathbb C}}
\renewcommand {\b}{\beta}
\newcommand {\Q}{\mathcal Q}
\newcommand{\dop}[1]{\frac{\partial}{\partial #1}}
\newcommand{\dopt}[2]{\frac{\partial #1}{\partial #2}}
\newcommand{\vardop}[3]{\frac{\partial^{|#3|} #1}{\partial {#2}^{#3}}}
\newcommand{\tc}[1]{T^{\mathbb C}_{#1}}
\newcommand{\ph}{\varphi}

\newcommand{\CC}[1]{\mathbb{C}^{#1}}
\newcommand{\CP}[1]{\mathbb{CP}^{#1}}
\newcommand{\RR}[1]{\mathbb{R}^{#1}}
\newcommand{\dw}{\frac{\partial}{\partial w}}
\newcommand{\dz}{\frac{\partial}{\partial z}}
\numberwithin{equation}{section}
\newcommand{\diffcr}[1]{\rm{Diff}_{CR}^{#1}}
\newcommand{\Hol}[1]{\rm{Hol}^{#1}}
\newcommand{\Aut}[1]{\rm{Aut}^{#1}}
\newcommand{\hol}[1]{\mathfrak{hol}^{#1}}
\newcommand{\aut}[1]{\mathfrak{aut}^{#1}}
\newcommand{\fps}[1]{\C\llbracket #1 \rrbracket}
\newcommand{\fpstwo}[2]{#1\llbracket #2 \rrbracket}
\newcommand{\cps}[1]{\C\{#1\}}

\title[]{Sphericity and Analyticity of a strictly pdeusoconvex hypersurface in low regularity I}

\author{Ilya Kossovskiy}
\address{(i) Department of Mathematics $\&$ International Center of Mathematics, Sustech University, Shenzhen, China; (ii) Department of Mathematics and Statistics, Masaryk University, Brno, Czech Republic; (iii) TU Wien, Vienna, Austria}
\email{ilyakos@sustech.edu.cn}
\author{Dmitri Zaitsev}
\address {School of Mathematics, Trinity College Dublin, Dublin, Ireland}
\email{zaitsev@maths.tcd.ie}





\begin{abstract}
In our earlier work \cite{KZ}, we introduced an {\em analytic regularizability theory} for smooth strictly pseudoconvex hypersurfaces in complex space. That is, we found a necessary and sufficient condition for a hypersurface to be CR-equivalent to an {\em analytic} target. The condition amount to the holomorphic extension property for a smooth function on a totally real submanifold, both the function and the submanifold being uniquely associated with the given hypersurface.

In the present paper, we develop our method further. First, we extend the result in \cite{KZ} to hypersurfaces of finite (possibly low) smoothness.  Second, we introduce a new tool for studying CR hypersurfaces in low regularity called {\em regularizing $(0,1)$ sections}. Using the latter key tool, we solve the open problem of checking the {\em sphericity} of a strictly pseudoconvex hypersurface in $\CC{2}$ in {low regularity}, precisely in the  regularity $C^k,\,2\leq k< 7$ which  is {\em not} covered by the classical Cartan-Tanaka-Chern-Moser theory. As an application of our theory, we deduce the sphericity of a strictly pseudoconvex hypersurface in $\CC{2}$ of regularity $C^6$ with vanishing Cartan-Chern CR-curvature.
\end{abstract}

\maketitle   

\def\Label#1{\label{#1}{\sf (#1)}~}


\def\cn{{\C^n}}
\def\cnn{{\C^{n'}}}
\def\ocn{\2{\C^n}}
\def\ocnn{\2{\C^{n'}}}


\def\dist{{\rm dist}}
\def\const{{\rm const}}
\def\rk{{\rm rank\,}}
\def\id{{\sf id}}
\def\aut{{\sf aut}}
\def\Aut{{\sf Aut}}
\def\CR{{\rm CR}}
\def\GL{{\sf GL}}
\def\Re{{\sf Re}\,}
\def\Im{{\sf Im}\,}
\def\span{\text{\rm span}}

\def\codim{{\rm codim}}
\def\crd{\dim_{{\rm CR}}}
\def\crc{{\rm codim_{CR}}}

\def\lr{\longrightarrow}
\def\phi{\varphi}
\def\eps{\varepsilon}
\def\d{\partial}
\def\a{\alpha}
\def\b{\beta}
\def\g{\gamma}
\def\G{\Gamma}
\def\D{\Delta}
\def\Om{\Omega}
\def\k{\kappa}
\def\l{\lambda}
\def\L{\Lambda}
\def\z{{\bar z}}
\def\w{{\bar w}}
\def\Z{{\1Z}}
\def\t{\tau}
\def\th{\theta}

\def\H{\hat H}

\emergencystretch15pt
\frenchspacing

\newtheorem{Thm}{Theorem}[section]
\newtheorem{Cor}[Thm]{Corollary}
\newtheorem{Pro}[Thm]{Proposition}
\newtheorem{Lem}[Thm]{Lemma}

\theoremstyle{definition}\newtheorem{Def}[Thm]{Definition}

\theoremstyle{remark}
\newtheorem{Rem}[Thm]{Remark}
\newtheorem{Exa}[Thm]{Example}
\newtheorem{Exs}[Thm]{Examples}

\def\bl{\begin{Lem}}
\def\bl{\begin{Lem}}
\def\el{\end{Lem}}
\def\bp{\begin{Pro}}
\def\ep{\end{Pro}}
\def\bt{\begin{Thm}}
\def\et{\end{Thm}}
\def\bc{\begin{Cor}}
\def\ec{\end{Cor}}
\def\bd{\begin{Def}}
\def\ed{\end{Def}}
\def\br{\begin{remark}}
\def\er{\end{remark}}
\def\be{\begin{Exa}}
\def\ee{\end{Exa}}
\def\bpf{\begin{proof}}
\def\epf{\end{proof}}
\def\ben{\begin{enumerate}}
\def\een{\end{enumerate}}
\def\beq{\begin{equation}}
\def\eeq{\end{equation}}

\section{Introduction}

Let $M\subset\CC{n+1},\,n\geq 1$, 
be a strictly pseudoconvex real hypersurface.
We call $M$ {\em spherical} at a point $p\in M$
 if there exists a
CR-diffeomorphism 
between neighborhoods of $p$ in $M$
and of a point $p'$ in the standard sphere $S^{2n+1}\subset\CC{n+1}$. 
The celebrated Cartan-Tanaka-Chern-Moser theory \cite{cartan,tanaka,chern} provides a compete answer to the problem of shericity of a 
{\em smooth}
strictly pseudoconvex hypersurface. This is done by building the so-called {\em canonical Cartan connection} associated with $M$  and its CR structure, 
and subsequently computing its curvature. 
The sphericity at a point $p$ is then 
equivalent to the vanishing
in a neighborhood of $p$ of certain principal curvature components
of the Cartan connection.
 However, 
 already
 the 
 construction of the Cartan connection on $M$ and its curvature
 requires the regularity of at least $C^7$ (as can be verified by the construction in each of the cited paper, and also those relying on the theory of {\em Parabolic Geometries}, see e.g. \cite{cs}). 

 We shall note here that, even though the computation of the principal curvature components in the Cartan-Tanaka-Chern-Moser theory requires derivatives of the defining function of respectively the order $\leq 6$ ($n=1$) or the order $\leq 4$ ($n\geq 2$), concluding the sphericity of a hypersurfaces with vanishing principal curvature components requires the {\em full connection and its curvature} (which relies on derivatives of the defining function of order $\leq 7$, as discussed above).
 
%
%
In this way, the corresponding sphericity problem in low regularity 
remained  open since the 1970's:

\medskip

\noindent {\bf Open Problem.} 
Find necessary and sufficient conditions for 
a strictly pseudoconvex real hypersurface
$M\subset\CC{n+1}$
 of class $C^k$, where $k<7$,
to be spherical at a point $p\in M$.
\medskip 

More precisely, we shall adopt the following terminology:

\bd
A real hypersurface $M\subset\CC{n+1}$
of class $C^k,\,k\geq 2$
is called {\em spherical of class $C^s$}
at a point $p\in M$, 
$s\le k$,
if there exists a $C^s$ CR-diffeomorphism
between a neighborhood of $p$ in $M$
and an open subset of the standard sphere
$S^{2n+1}\subset\CC{n+1}$.
We call $M$ {\em spherical at $p$}
(without specifying the regularity class
of the CR-diffeomorphism)
if it is spherical at $p$ of class $C^0$.
\ed
As discussed in Section 2, one can assume $s$ here to be as large as $s=k-1/2$ (unless $k-1/2\in\mathbb N$, in which case $s$ can be chosen as $s=k-1/2-\epsilon$ with any small $\epsilon>0$).

The above Open Problem is closely related to numerous theorems in Complex Analysis asserting that domains with spherical boundaries are biholomorphic to a ball or are universally covered by a ball (and then are certain ball quotients), see e.g. Chern-Ji \cite{chernji}, Huang-Ji \cite{huangji}, Nemirovski-Shafikov \cite{ns,ns2}, Li-Luk \cite{liluk}, Huang-Xiao \cite{huangxiao}. 
These results often assume 
merely $C^2$ regularity of the boundary, 
while {\em no} criterion for inspecting the sphericity
based on the Cartan connection
 is available in such a low regularity!
 
 More broadly, understanding low regularity
 in problems of geometry and PDEs
 is known to be difficult
 and brings new insights.

The main goal of the present paper is to provide 
a solution to the above Problem
for 
$ M\subset\CC{n+1}$ of class $C^k$ with $ k\geq 2$ in the dimension $n+1=2$,
see \autoref{main-2} below. Higher dimensions $n+1\geq 3$ are to be addressed in the upcoming paper \cite{KZnext} (which is mainly due to the fact that the Cartan-Tanaka-Chern-Moser theory looks quite different respectively in dimension $2$ and in higher dimensions). 
To obtain sharper regularity results,
we also consider H\"older classes $C^k$ 
for arbitrary nonnegative number $k\in\R$.
Our sphericity criterion amounts to an (explicit) system of
 differential equations for the defining function of a hypersurface,
  and the requirement of the higher than expected regularity of 
  several functions uniquely associated with the hypersurface 
  (the latter can be computed explicitly in terms of the defining function).    

Our approach relies on
the recent {\em analytic regularizability theory} 
by the authors in \cite{KZ}. 
That theory provides a necessary and sufficient condition for a
$C^\infty$ {\em smooth}
 strictly pseudoconvex hypersurface ({\em Condition E}, where ``E'' stands for ``extension'') to be (locally) 
CR-diffeomorphic to a real-analytic hypersurface.
%
Condition E is formulated in terms of a certain 
{\em holomorphic extension property}
 for a distinguished function invariantly associated with the hypersurface (called here the {\em $\Phi$-function}). 
 This condition 
 established in \cite{KZ} in the $C^\infty$ category, 
 is here extended to the case of low regularity (see Section 2) and subsequently used for addressing the sphericity problem.
 An important role is played by the new 
tool introduced in this paper that we call
 {\em regularizing vector fields}. 
 The latter are certain $(0,1)$ vector fields
 whose differentiation 
 of the $1$-jet section 
 induced by the complex tangent bundle 
 $T^\C M$
 gives CR functions of higher than expected regularity. 
%

More precisely, 
consider the holomorphic fiber bundle $J^{1,n}$
of complex hyperplanes over $\CC{n+1}$
and its canonical section over $M$
\beq
\Label{th-section}
s_M^1\colon M\to J^{1,n}, 
\quad
q\mapsto (q, [T^\C_q M]).
\eeq
Fixing a system of local holomorphic  coordinates
in $\CC{n}\times\CC{}$ yields a (local) trivialization of the bundle $J^{1,n}$, and we can write
\beq
\Label{th-def}
s_M^1(q) = (q, \th(q))
\in M\times \CC{n} \subset M\times \4P^n,
\eeq
where 
we use the standard embedding of $\C^n$ into the projective space and
$\th=(\th_1,\ldots,\th_n)$ is of class $C^{k-1}$.
It is easy to express $\th$
in terms of a local defining function $\rho$ of $M$
with $\rho_w\ne 0$
in local holomorphic coordinates
\beq
\Label{th-coor}
\th_j = -\frac{\rho_{z_j}}{\rho_w}, 
\quad j=1,\ldots,n,
\quad
(z,w)=(z_1,\ldots,z_n,w)\in \CC{n}\times\C.
\eeq

We now introduce our first new tool
for dealing with low regularity.
When writing $C^k$,
we shall always consider arbitrary real $k\ge0$, and all vector fields are assumed to be of class $C^0$
unless specified otherwise.

\bd
\Label{regu}
Let $M\subset\CC{n+1}$
be a real hypersurface of class $C^k$, $k\ge 2$,
with $\th$ as in \eqref{th-def} or \eqref{th-coor}.
A $(0,1)$ vector field $L$ on $M$
is called {\em regularizing of class $C^s$}
(or $C^s$-regularizing), $s\le k$,
if $L\th$ is a (vector-valued) 
CR function on $M$ of class $C^s$.
\ed

Note that there is no assumption
on the regularity of $L$ itself,
only on the derivative $L\th$ of $\th$ along $L$
(in fact, certain regularity follows automatically,
see \autoref{basis-l} below).
%
%
Also note
that regularizing vector fields $L$
are defined for any class $C^s$ with $s\le k$,
even if $L$ itself could be of lower regularity.

\br
\Label{reg-basis}
It is easy to construct
regularizing $(0,1)$ vector fields
on a
Levi-nondegenerate pseudoconvex hypersurface $M$.
Indeed, choose any basis $L_1,\ldots,L_n$
of $(0,1)$ vector fields 
of class $C^{k-1}$
and normalize
them to satisfy 
\beq\Label{kron}
L_j \th_l = \delta_{jl}
\eeq
(where we use the Kronecker $\delta$). This is always possible
by multiplying with the inverse of the 
matrix $(L_j\th_k)$ 
(which is of class $C^{k-2}$), 
the latter matrix
being invertible
 in view of the Levi-nondegeneracy.
 Since the constant functions $\delta_{jl}$
are CR and of class $C^s$ on $M$ for any $s\le k$,
all $L_j$ are $C^s$-regularizing.
So constructed $C^s$-regularizing
vector fields are themselves of class $C^{k-2}$,
which is allowed to be lower than $s$.
\er

\begin{definition}
The particular system of regularizing vector fields constructed in \autoref{reg-basis} is called {\em (the system of) basic regularizing $(0,1)$ sections}.    
\end{definition}

\br
\Label{basis-l}
Given any basis $L_1,\ldots,L_n$
of regularizing $(0,1)$ vector fields
(that always exists in view of Remark~\ref{reg-basis})
for a Levi-nondegenerate hypersurface $M$,
an arbitrary $(0,1)$ vector field
$L$ is regularizing of class $C^s$
precisely whenever its expansion
$L=c_1 L_1+\ldots+c_n L_n$
has all coefficients $c_j$
being CR functions of class $C^s$.
Indeed, since
$$
\begin{pmatrix}
L\th_1\\
\vdots\\
L\th_n\\
\end{pmatrix}
=
\begin{pmatrix}
L_1\th_1 & \cdots & L_n \th_1\\
\vdots & \ddots & \vdots\\
L_1\th_n & \cdots & L_n \th_n\\
\end{pmatrix}
\begin{pmatrix}
c_1\\
\vdots\\
c_n\\
\end{pmatrix}
$$
and the matrix
$(L_j\th_l)$
is CR of class $C^s$
and
 is invertible 
by the Levi-nondegeneracy,
$L\th$ is CR of class $C^s$
if and only if the coefficients 
$c_j$ are all  CR of class $C^s$.
In particular, using the basis
of $C^s$-regularizing $(0,1)$
vector fields 
of class $C^{k-2}$
constructed in 
\autoref{reg-basis},
it follows that any $C^s$-regularizing $(0,1)$
vector field $L$ is automatically
of class (at least) $C^{\min(s,k-2)}$.
\er

To formulate our sphericity criteria in low regularity,
we also need the $2$-jet version of the above $\th$-function,
i.e.\
the function $\Phi$ valued in symmetric $n\times n$ matrices
such that
\beq
\Label{Phi-def}
j^2_q Q_q = (q, \th(q), \Phi(q)) 
\in M\times \C^n \times \C^{\frac{n(n+1)}2}
\subset
J^{2,n},
\eeq
where $J^{2,n}$ is the holomorphic fiber bundle
of $2$-jets of complex hypersurfaces in $\CC{n+1}$
and
$j^2_q Q_q$ is the $2$-jet of the formal Segre variety
that is well-defined for $k\ge 2$.
Similarly to \eqref{th-coor}, 
the function $\Phi$ can be calculated
in terms of a local defining function $\rho$ of $M$
with $\rho_w\ne 0$
by the formula
\begin{equation}\Label{Phi00}
\Phi = (\Phi_{ij}),
\quad
	\Phi_{ij}=\frac{1}{(\rho_w)^3}
\begin{vmatrix} 0 & \rho_{z_j} & \rho_w\\ \rho_{z_i} & \rho_{z_iz_j} & \rho_{z_iw} \\ \rho_w & \rho_{z_jw} & \rho_{ww} \end{vmatrix},
 \quad i,j=1,...,n.
\end{equation}

We are now ready to formulate
our first sphericity result in $\C^2$. 

\begin{theorem}\Label{main-2}
Let
$M\subset\CC{2}$
be  a strictly pseudoconvex  real hypersurface 
of class $C^k,\,k\geq 2,$
and $L$ a (nonzero) regularizing
$(0,1)$
vector field on $M$ of class $C^s$.
Then the following holds:

\smallskip

(i) if $2\leq s\leq k$ and $M$ is spherical of class $C^s$, then one has
\beq
\Label{main-cond}
\Phi,L\Phi, L^2\Phi, L^3\Phi \in C^{s-2},
\quad
L^4 \Phi =0.
\eeq

\smallskip

(ii) if $1\leq s \leq k$ and the conditions 
\beq
\Label{main-cond'}
\theta,\Phi,L\Phi, L^2\Phi, L^3\Phi \in C^{s},
\quad
L^4 \Phi =0.
\eeq
are satisfied, then $M$ is spherical of class $C^{(s+1)-}$.
\end{theorem}






\autoref{main-2} has the following two immediate corollaries.  

First, for the regularity $k> 7/2$, we are able to formulate a {\em complete criterion} of sphericity:

\begin{theorem}
\Label{main-2-iff}
Let
$M\subset\CC{2}$
be a strictly pseudoconvex real hypersurface 
of class $C^k,\,k> 7/2$
and $L$ a (nonzero) regularizing
$(0,1)$
vector field on $M$ of class $C^1$.
Then $M$
is spherical
at a point $p$
 if and only if
\beq
\Label{main-cond1}
L\Phi, L^2\Phi, L^3\Phi \in C^1,
\quad
L^4 \Phi =0.
\eeq
\end{theorem}

Second, in the "super low" regularity $2\leq k\leq 7/2$, we can formulate a sufficient condition for sphericity.

\begin{corollary}
\Label{main-2-suff}
Let
$M\subset\CC{2}$
be a strictly pseudoconvex real hypersurface 
of class $C^k,\,2\leq k\leq 7/2,$
and $L$ a (nonzero) regularizing
$(0,1)$
vector field on $M$ of class $C^1$.
Then $M$
is spherical
at a point $p$
 if 
\begin{equation}
\Label{main-cond2}
\Phi,L\Phi, L^2\Phi, L^3\Phi \in C^1,
\quad
L^4 \Phi =0.
\end{equation}
\end{corollary}

\br
\Label{cond-indep}
Another important feature of 
our conditions in \autoref{main-2}
is their independence of the choice of $L$.
Indeed, since $L$ is nonzero,
it forms a basis of $(0,1)$ vector fields
(in $\C^2$). It can be verified from Remark ~\ref{basis-l}: in the $\CC{2}$ case,
Remark~\ref{basis-l}
implies that any other 
regularizing vector field $\2L$
of class $C^s$
is a multiple of $L$
by a CR function of class $C^s$
and since $(0,1)$ vector fields
commute with CR functions,
condition \eqref{main-cond}
remains invariant
when $L$ is replaced by $\2L$. 
\er


Finally, we have an application of \autoref{main-2} directly to the Cartan-Tanaka-Chern-Moser theory (in low regularity). Precisely, we prove the following result with optimal regularity:

\begin{theorem}\Label{zerocurv}
A strictly pseudoconvex hypersurface $M\subset\CC{2}$ of class $C^6$ is spherical if and only if its Cartan invariant (the CR curvature) vanishes on $M$.    
\end{theorem}

Note that $6$ is precisely the number of derivatives needed to define the CR curvature in $\C^2$, 
hence the regularity order $6$ is optimal.
On the other hand,
as discussed above, this result can {\em not} be deduced from Cartan-Tanaka-Chern-Moser theory,
which requires higher order differentiation
to conclude sphericity by standard methods.

\bigskip

\begin{center} \large \bf Acknowledgements          \end{center}

\mbox{}

The first  author was supported by the GACR grant 22-15012J, the FWF grant P34369, and the NSFC grant W2431006.

\bigskip

\section{Condition E for smoothly embedded real hypersurfaces}

 \subsection{The method of associated differential equations}
 The study of the relationship between the geometry of 
 real hypersurfaces in $\C^2$ and 2nd order ODEs 
  \begin{equation}\Label{wzz}
 w\rq{}\rq{}=\Phi(z,w,w\rq{}).
 \end{equation}
 goes back to Segre \cite{segre}
 and Cartan \cite{cartan},
 see also Webster \cite{webster,We00}
More generally,  the geometry of a real hypersurface in $\CC{n+1},\,n\geq 1$, is related to
 that of a complete second order system of PDEs
\begin{equation}\Label{wzkzl}
w_{z_kz_l}=\Phi_{kl}(z_1,...,z_n,w,w_{z_1},...,w_{z_n}),\quad \Phi_{kl}=\Phi_{lk},\quad k,l=1,...,n,
\end{equation}
  Moreover, {\em in the real-analytic case},
  this relation becomes 
more
explicit 
by means of the Segre family.
Namely,
to any real-analytic Levi-nondegenerate hypersurface $M\subset\CC{n+1},\,n\geq 1$, one can uniquely associate a holomorphic ODE \eqref{wzz} ($n=1$) or a holomorphic PDE system \eqref{wzkzl} ($n\geq 2$),
whose solutions are precisely the Segre varieties.
The Segre family of $M$ plays a role of a 
``mediator'' between the hypersurface and the associated differential equations.  
For recent work on associated differential equations in the degenerate setting, see e.g. the papers  
  \cite{divergence, nonminimalODE, nonanalytic} of the first author with Lamel and Shafikov.

  The associated differential equation procedure is particularly simple in the case of a Levi-nondegenerate hypersurface in $\CC{2}$. In this case the Segre family is an
 anti-holomorphic 
 2-parameter family of complex holomorphic curves. 
 It then follows from the standard ODE theory that there exists a unique ODE \eqref{wzz}, for which the Segre varieties are precisely the graphs of solutions. This ODE is called \it the associated ODE. \rm

In general, both right hand sides in \eqref{wzz},\eqref{wzkzl} appear as functions determining the $2$-jet of a Segre variety 
at a given point
as a function of the $1$-jet at the same point.
More explicitly, we use coordinates
$$
	(z,w)=(z_1, \ldots, z_n,w) \in \CC{n}\times \C = \CC{n+1}.
$$ 
Fix a 
 real-analytic
hypersurface
 $M\subset\CC{n+1}$
 passing through the origin, 
 and choose a 
sufficiently small neighborhood $U$
 of the origin.
  When $M$ is Levi-nodegenerate, we can associate a 
{\em complete second order system of holomorphic PDEs}
 to $M$,
which is uniquely determined by the condition that the differential equations are satisfied by all the
graphing functions $h(z,\zeta) = w(z)$ of the
Segre varities $\{Q_\zeta\}_{\zeta\in U}$ of $M$ in a
neighbourhood of the origin.

To be more precise, we consider a
so-called {\em  complex defining
 equation } (see, e.g., \cite{ber}),
$$
w=\rho(z,\bar z,\bar w),
$$ 
of $M$ near the origin, which one
obtains by substituting 
$u=\frac{1}{2}(w+\bar w),\,v=\frac{1}{2i}(w-\bar w)$ 
into 
a real-analytic defining equation and
solving for $w$ as function of $(z,\bar z, \bar w)$
by the implicit function theorem.
 The Segre
variety $Q_x$ of a point 
$$x=(a,b)\in U,\quad
a\in\CC{n},\,b\in\CC{},$$ 
is  now given
as the graph of the function
\begin{equation} 
\Label{segredf}
w (z)=\rho(z,\bar a,\bar b), 
\end{equation}
where we slightly abuse the notation
using the same letter $w$ for 
both the last coordinate and
the function $w(z)$ defining a Segre variety.
Differentiating \eqref{segredf} we obtain
\begin{equation}\Label{segreder} 
	w_{z_j}=\rho_{z_j}(z,\bar a,\bar b),
	\quad
	j=1,\ldots,n. 
\end{equation}
Considering \eqref{segredf} and \eqref{segreder}  as a holomorphic
system of equations with the unknowns $\bar a,\bar b$, 
in view of the Levi-nondegeneracy of $M$,
an
application of the implicit function theorem yields holomorphic functions
 $A_1,...,A_n, B$ such that
 \eqref{segredf} and \eqref{segreder} are solved by
$$
	\bar a_j=A_j(z,w,w'),\quad
	\bar b=B(z,w,w'),
$$
where we write
$$
	w' = (w_{z_1},  \ldots, w_{z_n}).
$$
The implicit function theorem applies here because the
Jacobian of the system coincides with the Levi determinant of $M$
for $(z,w)\in M$ (\cite{ber}). Differentiating \eqref{segredf} twice
and substituting the above solution for $\bar a,\bar b$ finally
yields
\begin{equation}\Label{segreder2}
w_{z_kz_l}=\rho_{z_kz_l}(z,A(z,w,w'),
B(z,w,w'))=:\Phi_{kl}(z,w,w'),
\quad
k,l=1, \ldots, n,
\end{equation}
or, more invariantly,
i.e.\ independent of the coordinate choice,
\begin{equation}\Label{segreder2'}
	j^2_{(z,w)} Q_x = \Phi(x, j^1_{(z,w)} Q_x).
\end{equation}
Now \eqref{segreder2}
(or \eqref{segreder2'})
is the desired complete system of holomorphic second order PDEs
denoted by $\mathcal E = \mathcal{E}(M)$.

 \begin{definition}\Label{PDEdef}
 We call $\mathcal E = \mathcal{E}(M)$  \it the system of PDEs 
 associated with $M$. \rm  
\end{definition}

\subsection{Condition E}
Our main results rely on a 
new
{\em analytic regularizability} condition
for hypersurfaces of low regularity.
In the $C^\infty$ case,
such condition was
 introduced in \cite{KZ} and  called    
  {\em Condition E (``E'' for extension)}.
 We shall give below both an invariant 
 and a coordinate-based formulations of it. 
 For the basic concepts in CR Geometry (such as Segre varieties and formal submanifolds) we refer to \cite{ber}, and for {\em jet bundles} and related concepts to \cite{cs}. 

Let 
$$
	\pi \colon J^{1,n}
	\to \CC{n+1},
$$ 
be the bundle
of $1$-jets of complex hypersurfaces of $\CC{n+1}$, 
which is a projective holomorphic 
 fiber bundle over $\CC{n+1}$ with the fiber 
isomorphic to $\mathbb P^n$,
and let
$M\subset\CC{n+1},\,n\geq 1$, 
be (for the moment) a ($C^\infty$) smooth
{\em  strictly pseudoconvex} real hypersurface. 
Then the complex tangent bundle $\tc{}M$ induces
the natural smooth (global) embedding 
$$
	\varphi\colon M\to J^{1,n}, \quad x \mapsto \bigl(x,[\tc{x}M]\bigr).
$$ 
The image 
$$
	\ph(M)=:M_J \subset J^{1,n}
$$ 
 is a  smooth totally real
 $(2n+1)$-dimensional real submanifold  in the $(2n+1)$-dimensional complex manifold $J^{1,n}$
 by an observation of
 Webster \cite{webster}.
Next, 
associated with $M$ is the  smooth (weakly) pseudoconvex real hypersurface
$$
	\widehat M:=\pi^{-1}(M) \subset J^{1,n}.
$$
The manifold $M_J$ is a smooth real submanifold in $\widehat M$.  
Note that  $\widehat M$ itself is locally CR-equivalent 
to $M\times \CC{n}$ (and thus is {\em holomorphically degenerate}, see \cite{ber}).  
In what follows we denote by $U^+$ 
the {\em pseudoconvex side of $M$},
by which we mean
the subset of an open neighborhood
$U$
of a point of $M$ 
defined by $\rho<0$, i.e.\
$$
U^+ = U \cap \{\rho<0\},
$$
where $\rho$ is a local defining function of $M$ in $U$
with $d\rho\ne0$ and the
 complex hessian satisfying
$\d\bar\d\rho(X,\1X)\ge0$
for $X\in T^{10}M$.
Given such $U^+$, 
we write
 $$\widehat U^+:=\pi^{-1}(U^+).$$ 
Then $\3M$
is obviously (weakly) pseudoconvex
and $\3U^+$ is the pseudoconvex side of $\3M$.
Since all our considerations are local,
the exact choice of neighborhoods
of the reference point won't play any role.

We next fix a point $p\in M$
along with the corresponding point
$$
\3p:= (\pi|_{M_J})^{-1}(p) \in M_J \subset \3M\subset J^{1,n}.
$$
Since $M$ is smooth, we may consider at each point $q\in M$ 
near $p$, its formal complexfication at $q$ as a formal complex hypersurface in $\CC{n+1}\times\overline{\CC{n+1}}$ obtained by complexifying the formal Taylor series
at $q$ of the defining function $\rho$.
In this way, the formal Segre variety $Q_q$ of $M$ 
is defined as a formal complex hypersurface
defined by power series in $\CC{}[[Z-q]]$,
where $Z=(z,w)\in \CC{n}\times\C$ as before.
Then the map
$2$-jets 
\begin{equation}\Label{2jets}
q\in M \mapsto j^2_qQ_q \in J^{2,n}
\end{equation} 
induces a {\em smooth embedding}
of $M$ 
 into the bundle  
$$J^{2,n}=J^{2,n}(\CC{n+1})$$ 
of $2$-jets of complex hypersurfaces in $\CC{n+1}$. The space $J^{2,n}$  is canonically a fiber bundle  
$$
\pi^2_1 \colon J^{2,n}\to J^{1,n}
$$
over the $1$-jet bundle $J^{1,n}$.
The $2$-jet embedding
\eqref{2jets}
 defines a canonical section of $\pi^2_1$, 
\beq
\Label{s12}
s_1^2\colon M_J\to J^{2,n}.
\eeq

Now we recall
 our analytic regularizability condition for a  smooth strictly pseudoconvex hypersurface \cite{KZ}:

\begin{definition}\Label{coorfree} 
We say that $M$ {\em satisfies Condition E at $p\in M$}
 if
 the canonical section $s_1^2$
 given by \eqref{s12} extends 
 holomorphically and smooth up to the boundary
 to a neighborhood of $\3p$
 in
  the pseudoconvex side 
 $\widehat U^+ \cup \widehat M$.
\end{definition}


We next give an (equivalent to the above) coodinate-based formulation of Condition E. 
Let $M\subset\CC{n+1}$ be a smooth hypersurface with the defining equation 
\begin{equation}\Label{rhoeq}
\rho(Z,\bar Z)=0,\quad Z=(z,w)=(z_1,...,z_n,w)\in\CC{n+1},
\end{equation}
and
$p\in M$ a fixed point where
$\rho_w(p,\bar p)\neq 0$. 
Then the formal Segre variety at a point $q=(\tilde q,q_{n+1})\in M$ near $p$ is a graph of a function $w(z)$ (considered as a formal power series in $(z-\tilde q)$). Then the $2$-jets \eqref{2jets} are represented by either the scalar function 
$\Phi$ defined pointwise as $w''(z)$ for $z=\tilde q$  (case $n=1$),  or the symmetric matrix function $\Phi=(\Phi_{ij}),\,\,i,j=1,...,n$, defined pointwise as the collection of the partial derivatives
 $w_{z_iz_j}$ for $z=\tilde q$ (the general case $n\ge1$). It is possible to verify that, in turn, for $n=1$ we have
\begin{equation}\Label{Phi1}
	\Phi=\frac{1}{(\rho_w)^3}
	\begin{vmatrix} \rho & \rho_z & \rho_w\\ \rho_z & \rho_{zz} 
	& \rho_{zw} \\ \rho_w & \rho_{zw} & \rho_{ww} \end{vmatrix},
\end{equation}
and for $n>1$ we have 
\begin{equation}\Label{Phi}
	\Phi_{ij}=\frac{1}{(\rho_w)^3}
\begin{vmatrix} \rho & \rho_{z_j} & \rho_w\\ \rho_{z_i} & \rho_{z_iz_j} & \rho_{z_iw} \\ \rho_w & \rho_{z_jw} & \rho_{ww} \end{vmatrix}, \quad i,j=1,...,n.
\end{equation}
 (To obtain \eqref{Phi1},\eqref{Phi}, one has to differentiate the identity \eqref{rhoeq} once, assuming $w$ to be a function of $z$, express all the  $w_{z_j}$ in terms of the $1$-jet of $\rho$,
 and then differentiate \eqref{rhoeq} once more to obtain $w_{z_iz_j}=\Phi_{ij}$ in terms of the $2$-jet of $\rho$.)
  Both the scalar function \eqref{Phi1} and the matrix valued function \eqref{Phi} can be considered as either smooth functions on the strictly pseudoconvex hypersurface $M$ or as that on the totally real submanifold $M_J\subset J^{1,n}$ introduced above.

In terms of the $\Phi$-function \eqref{Phi1}-\eqref{Phi}, 
Condition E reads in an equivalent form as follows:
\begin{definition}\Label{coor} 
We say that $M$ {\em satisfies Condition E at $p$}, if 
 the function $\Phi$ defined on $M_J$ by either \eqref{Phi1} or \eqref{Phi} extends  
   holomorphically and smoothly up to the boundary
   to a neighborhood of $\3p$ in
 the pseudoconvex side $\widehat U^+ \cup \widehat M$.
\end{definition}
It is obvious that \autoref{coor} is equivalent to \autoref{coorfree}. 
Now the main result of \cite{KZ} is formulated as follows:

\begin{theorem}
\Label{mainold} 
A  smooth strictly pseudoconvex  real hypersurface $M\subset\CC{n+1},\,n\geq 1$,
 is {\em  analytically regularizable
at a point $p\in M$},
i.e.\ 
 a neighborhod of $p$ in $M$
 is smoothly CR-diffeomorphic
  to a real-analytic  hypersurface in
$\CC{n+1}$,
 if and only $M$ satisfies  Condition E at $p$.
\end{theorem}

Our first goal is to generalize \autoref{mainold} 
to $M$ of low regularity. 
We shall first formulate Condition E in that setting. 
As customary, for a real $k\ge0$,
we consider the class $C^k$
of functions 
whose derivatives exist and are continuous
up to order $[k]$ 
(where $[k]$ denotes the integral part of $k$),
and whose
derivative up to order $[k]$
are either continuous for $k=[k]$
or
in the local H\"older class $C^{k-[k]}$
for non-integer $k$.

\begin{definition}
We say that a strictly pseudoconvex hypersurface
 $M\subset\CC{n+1}$ of class $C^k,\,k\geq 2$,
 is {\em $C^s$ analytically regularizable at $p\in M$}
 for $k\geq s\geq 0$, if there exists a $C^s$ CR-diffeomorphism of a neighborhood of $p$ in $M$ into an analytic strictly pseudoconvex hypersurface.
\end{definition}

\begin{definition}
Given a function class $\6S$,
we say that a strictly pseudoconvex hypersurface
 $M\subset\CC{n+1}$ of class $C^k,\,k\geq 2$,
 satisfies {\em Condition E of class $\6S$}
 at $p\in M$,
   if  for some choice of a neighborhood $U$ of $p$, the associated function $\Phi$ defined on $M_J$ by
   $2$-jets of Segre varieties,
   i.e.\
    either \eqref{Phi1} or \eqref{Phi},
     extends  to 
     a neighborhood of $\2p$
     of the pseudoconvex side $\widehat U^+ \cup \widehat M$ holomorphically and
  of class $\6S$ up to the boundary.
\end{definition}

\begin{remark}
\Label{huru}
Turning to regularizability, 
 a proper holomorphic map
 between 
 strictly pseudoconvex
 bounded domains  $D,D'\subset\CC{n+1}$
 with boundaries of class of 
 $C^k$, $k\ge2$,
 is known to have extension
to the closure $\1D$ of the class
 \beq
 \Label{rk}
 \begin{cases}
 C^{k-1/2} & k-1/2\notin\4Z
 \\
 \cap_{\eps>0} C^{k-1/2-\eps} &  k-1/2\in\4Z
 \end{cases}
 \eeq
  (see 
  Khurumov \cite{hurumov}, Coupet \cite{coupet}, and the survey 
  by Forstneric  \cite{forstneric}  in this regard). 
  This result is optimal even if assuming the target analytic. 
%
\end{remark}

To include both cases in \eqref{rk},
it will be convenient to use the notation
for the regularity class
$$
C^{a-} := \bigcup_{b\le a, \; b\notin\4Z} C^b.
$$
In particular, 
$C^{a-}=C^a$
whenever $a$ is not an integer.
Then 
the regularity class in \eqref{rk}
equals $C^{(k-1/2)-}$ in both cases.

Then
in our local situation of
a CR homeomorphism 
$H$
between 
one-sided neighborhoods of
points $p$ in $M$ and $p'\in M'$
with $H(p)=p'$,
where $M,M'$ are 
strictly pseudoconvex hypersurfaces
of class $C^k$,
a standard construction
of strictly pseudoconvex
domains $D,D'$
sharing open parts of their boundaries
with respectively
neighborhoods of $p$ in $M$
and of $p'$ in $M'$,
along with the above result,
shows that $H$
is of class $C^{(k-1/2)-}$.

%

We can now formulate our 
new analytic regularizability result in low
regularity.

\begin{theorem}\Label{maink} 
Let $M\subset\CC{n+1},\,n\geq 1$, 
be a  strictly pseudoconvex 
 real hypersurface  of regularity class 
 $C^k,\,k\geq 2$,
 and $p\in M$ is fixed. 
 Then the following holds:

%
\ben
\item[(i)] 
 if 
   $M$ 
 is analytically regularizable of class $C^s$   at $p$, where $2\leq s \leq k$,
 then $M$ satisfies  
 Condition 
  E 
   at $p$
   of class $C^{s-2}$.
%
%
\item[(ii)]
 if $M$ 
satisfies
Condition E at $p$ of class $C^s$
 with $1\leq s\leq k$,
 then $M$ is analytically regularizable  of class $C^{(s+1)-}$ at $p$.
 \een
\end{theorem}

Combining \autoref{maink} with the above mentioned regularity result for CR-diffeomorphisms
in Remark~\ref{huru} allows for 

\begin{corollary}
\Label{Maink} 
Let $M\subset\CC{n+1},\,n\geq 1$, 
be a  strictly pseudoconvex 
 real hypersurface  of regularity class $C^k,\,k\geq 2.$ 
Then the following holds:
\ben
\item
 if $k>5/2$ and a neighborhood of $p$ in $M$ 
 is CR-homeomorphic 
 to a strictly pseudoconvex
 real-analytic hypersurface,
 then it satisfies  
 Condition 
  E 
   at $p$
   of class $C^{(k-5/2)-}$.
 \item
 if $M$ satisfies
Condition E of class $C^1$
 at $p$,
 then  there exists a CR diffeomorphism
 of class $C^{(k-1/2)-}$
 from a neighborhood of $p$ in $M$
 onto a {\em strictly pseudoconvex}
 real-analytic
 hypersurface in $\CC{n+1}$.
\een
In particular, a neighborhood of $p$ in
a strictly pseudoconvex  real hypersurface 
$M\subset\CC{n+1}$
of class $C^k,\,k>  7/2$, 
is 
CR-homeomorphic to a
strictly pseudoconvex
real-analytic hypersurface in $\CC{n+1}$
 if and only if $M$ satisfies  Condition E
 at $p$ of class $C^1$.
\end{corollary}

\bpf
(1) 
Let $H$ be a CR homeomorphism
between a neighborhood of $p$ in $M$
and a strictly pseudoconvex
 real-analytic hypersurface $\2M\subset \CC{n+1}$.
By the results stated in \autoref{huru},
$H$ is of class $C^{(k-1/2)-}$
and 
the desired conclusion follows
from Theorem~\ref{maink} (i).

To show
(2) apply first Theorem~\ref{maink} (ii),
which which has the same assumptions.
The application yields 
 a CR diffeomorphism
 of class $C^{\min(s,k-1)-}$
 from a neighborhood of $p$ in $M$
 onto a {\em strictly pseudoconvex}
 real-analytic
 hypersurface in $\CC{n+1}$.
Now, by the result in \autoref{huru},
the CR diffeomorphism is of
class  $C^{(k-1/2)-}$ as claimed.
\epf

\begin{proof}[Proof of \autoref{maink}]

\smallskip

\noindent {\it Proof of (i):}

\smallskip
Let $H$ be the CR diffeomorphism.
Consider the Lewy extension of $H$,
also denoted by $H$ by a slight abuse of notation, 
to the pseudoconvex side $U^+$ of
(a neighborhood of $p$ in) $M$,
 which automatically has
 the same boundary  regularity,
 see e.g. \cite[Theorem~7.5.1]{ber}
 and the regularity properties
 of the solutions to the Dirichlet problem,
 see e.g. \cite{Go}.
 We then consider the 
 lift $H^{(2)}$ of $H$ 
 between the bundles $J^{2,n}$ 
 above $M$ and $\2M$ and its holomorphic extension
 over 
 $$
 \Omega:=(\pi^2)^{-1}(U^+),
 \quad
 \pi^2\colon J^{2,n}(\CC{n+1})\to \CC{n+1},
 $$ 
 lying over the pseudoconvex side $U^+$ of $M$, where $\pi^2$ is the canonical projection in $J^{2,n}$. 
 The lift $H^{(2)}$
 involves the 2nd order derivatives,
 hence has the boundary regularity of class $C^{s-2}$ in $\Omega$. 
 Similarly, the lift $H^{(1)}$
 between the $1$-jet bundles $J^{1,n}$
 is of class $C^{s-1}$.

 Pulling back the differential equations \eqref{wzz} or \eqref{wzkzl} in $J^{2,n}$  associated with the real-analytic 
 hypersurface $\2M$,
 we see that the $\Phi$-function of $M$
 equals the composition 
 $$
 \Phi=(H^{(2)})^{-1}\circ \2\Phi \circ H^{(1)},
 $$
 where $\2\Phi$ is the $\Phi$-function of $\2M$.
 Since $\2M$ is real-analytic, $\2\Phi$ is also real-analytic,
 hence $\Phi$ is of class $C^{s-2}$ as desired.
%

\smallskip

\noindent {\it Proof of (ii):} 

\smallskip

The proof here is analogous to the proof of
the sufficiency of Condition E in \autoref{mainold} 
given in \cite[\S4]{KZ} and is essentially obtained by repeating the arguments, while paying attention to the regularity classes.
 We only briefly outline the differences with the $C^\infty$ case:

\smallskip

- The main proof strategy is to change the complex structure
on the pseudoconcave side of $M$
by 
by identifying 
with
 the leaf space of certain ``Segre foliation''
 constructed using the extended $\Phi$-function.
 Then consider the local complex
 coordinates given by the low regularity
  Newlander-Nirenberg theorem
  due to
  Hill-Taylor \cite{ht},
  where the image of $M$
  becomes the fixed point set
  of an antiholomorphic involution,
  implying the desired analyticity.
  The involution is constructed
  on wedges with edge $M_J$
  and then extended by means of 
  the Edge-of-the-Wedge theorem.


\smallskip

- The Webster embedding map of $M$ into $J^{1,n}$
sending $q\in M$ into $[T^\C_qM]$
 has the (lower than that of $M$) regularity class $C^{k-1}$,
 while the function $\Phi$ has the regularity class $C^s$ by the assumption.

\smallskip

- As in \cite[Step I]{KZ},
the extension of $\Phi$ (again denoted by $\Phi$ by a slight abuse of notations) to $\3U^+\cup \3M$
defines the holomorphic
distribution $D$
in $\3U^+\subset J^{1,n}$ of class $C^s$ up to the boundary by
\begin{equation}\Label{distrib}
\omega=(\omega_0,\omega_1,...,\omega_n)=0,
\end{equation}
where 
\begin{equation}\Label{omegak}
\omega_0:=dw-\sum_{1}^n\xi_jdz_j,
\quad
\omega_k:=d\xi_k-\sum_1^n\Phi_{kl}dz_l,\,\,k=1,...,n.
\end{equation} 
Here $(z_1,...,z_n,w)=(z,w)$ denote the coordinates in $\CC{n+1}=\CC{n}\times\CC{}$ as before, 
and $\xi=(\xi_1,...,\xi_n)$ are the respective ``jet''-variables corresponding to the derivatives $w_{z_1},...,w_{z_n}$ respectively
as constructed in \cite[\S4.1]{KZ}.
Furthermore, $M_J$ is precisely the locus of points in $\widehat M$, for which $D$ {\em is tangent} (and hence complex tangent) to $\widehat M$. This follows directly from the construction of $M_J$ and $\widehat M$.

The proof of formal integrability
of $D$ in \cite[Proposition~4.1]{KZ}
only uses first order derivatives of $\Phi$,
hence still goes through in our case.
The respective holomorphic foliation is
called ``the one-sided Segre foliation'' $\mathcal F$.

- As in \cite[Step II, \S4.2]{KZ},
the foliation
   $\mathcal F$ admits 
   an extension $\mathcal F'$
       of regularity class $C^s$
    to a neighborhood in $J^{1,n}$ of
   the point 
   $$
   \3p:=(\pi|_{M_J})^{-1}(p)
   $$ (that is, the extension of $\mathcal F$ is given by level sets of a $C^s$ function),
    such that the leaves have connected intersections 
    with $\3U^+$
    and are open subsets of
    leaves
    of the extended foliation. This follows from e.g. the recent result by Yao \cite{yao}. 

\smallskip

- As in \cite[Step III, \S4.2]{KZ},
endow the space of leaves of $\mathcal F$
with the structure of a 
$(2n+2)$-manifold with boundary $M_J$ in the natural (quotient) topology. 
Indeed,
choose local coordinates
$(z,w,\xi)$
 vanishing at $\3p$
such that
 the tangent space at $0$ to 
the leaf of $\mathcal F'$ through $0$ is $w=0,\,\xi_j=0$
and such that the defining function of $M$
has positive definite real hessian.
Then the $(2n+2)$-subspace given by $z=0$
has the property that all the leaves of the extended foliation
$\mathcal F'$ intersect it transversally at single points each
(after possibly shrinking the given neighborhood of $\3p$). 
Thus the space of leaves of $\mathcal F'$ can be identified with a domain in $\RR{2n+2}$. 
Accordingly, the space of leaves of $\mathcal F$ is an open connected subset
 of the latter domain,
  given by the condition of having nonempty intersection with the 
  domain $\widehat U^+$. 
This gives a structure of a smooth manifold
of class $C^s$ on the space of leaves of $\mathcal F$
 in its natural quotient topology.

To show that the space of leaves of  $\mathcal F$ has, furthermore, a structure of a smooth $(2n+2)$-manifold with boundary that can be identified with $M_J$,
perform a local diffeomorphism
of class $C^s$ near $0$
 with the identity differential at $0$
 straightening the foliation $\6F$
 and keeping positive definite real hessian
 of the defining function of $M$.
Then in the new coordinates 
$(x_1,\ldots,x_{4n+2})$, 
the leaves become ``horizontal" $2n$-planes 
given by $$x_j=c_j,\,j=2n+1,...,4n+2,$$ 
where $c_j$ are constants, and the hypersurface $\widehat M$ becomes
$$
	x_{4n+2}=\varphi(x_1, \ldots ,x_{4n+1})
$$ 
for a smooth function $\varphi$
of class $C^s$
which
has positive definite real hessian in
the variables $x_1,...,x_{2n}$.
Now the condition that a leaf is tangent to 
$\widehat M$ at some point becomes
$$
c_{4n+2}=\varphi(x_1,...,x_{2n},c_{2n+1},...,c_{4n+1}),
\quad
\varphi_{x_j}(x_1,...,x_{2n},c_{2n+1},...,c_{4n+1})=0,
\quad j=1,..,2n,
$$   
proving that the leaf space $\3U^+/\6F$ of $\6F$
is a real $(2n+2)$-manifold with boundary
(identified with) $M_J$
of class $C^{s}$.


\smallskip

- As in \cite[Step IV, \S4.2]{KZ},
we note that the holomorphic foliation
$\6F$ on $\3U^+$
induces a complex structure
on its leaf space,
which we shall change to its conjugate
to match the complex structure
coming from the pseudoconvex side of $M$
(related to the antiholomorphic dependence
of Segre varieties on their parameter).
Since $\6F$ is of class $C^s$,
the complex structure on the leaf space
is of class $C^s$ up to the boundary.

- As in \cite[Step V, \S4.2]{KZ},
the projection $J^{1,n}\to \CC{n+1}$
defines an embedding of the leaves of $\6F$ 
as closed complex
hypersurfaces
in the pseudoconvex side $U^+$,
$$\mathcal S^+:=\Bigl\{\pi(T)\Bigr\}_{T\in \mathcal F},$$
where we write $T\in \mathcal F$
to mean that $T$ is a leaf of $\6F$.
%

We further note that the leaf space
 $\3U^+/\6F=:\mathcal U$ 
 is endowed with a natural $(n+1)$-dimensional anti-holomorphic family of 
 $n$-dimensional complex submanifolds
 $Q_x$
 given by all leaves $T$
 with projections $\pi(T)$ passing through a fixed point
 $x\in U^+$.
We denote the resulting family of submanifolds
$Q_x\subset\mathcal U$ by $\mathcal S^-$, which is an anti-holomorphic family parameterized by $x\in U^+$. 
For completeness of the picture, we also call, for each $y\in\mathcal U$, the respective leaf $T=y\in \mathcal S^+$ its {\em Segre variety} and denote the latter one by $Q_p$. We then obtain the following familiar symmetry property:
$$x\in Q_y \Leftrightarrow y\in Q_x, \quad x\in U^+,\,\,y\in \mathcal U.$$ 

- As in \cite[Step VI, \S4.2]{KZ},
we obtain a smooth manifold $U$ decomposed as a union of two manifolds with boundary of  class $C^{s}$:
$$U=(U^-\cup M)\cup (U^+\cup M),$$
where both $U^-$ and $U^+$ are endowed with their individual complex structures (for $U^-$ this is the integrable structure induced from $\mathcal U$ and for $U^+$ this is the standard complex structure induced from $\CC{n+1}$). Moreover, 
it is clear from the previous steps
that each of the two structures
 admits a smooth extension to the boundary 
of class $C^{s}$
(the standard complex structure on $U^+$
is of course real-analytic).

To show that the two structures 
agree on $M$,
%
let's approximate with real-analytic hypersurfaces,
for which both complex structures are standard.
An approximation of order $k$
provides coincidence of 
both complex structures.

%
%

- Now we use the low regularity version of the Newlander-Nirenberg theorem due to Hill-Taylor \cite{ht} to map
the complex structure on $U$ 
onto the standard complex structure. 
In view of \cite[Remark~1.7]{ht},
their theorem applies 
whenever the complex structure $J$
is of class $C^r$ for some $r>1/2$,
i.e.
$
J \in \cup_{r>1/2} C^r,
$
which holds in our case.
Furthermore,  by \cite[Proposition~5.1]{ht},
the map providing holomorphic coordinates
is, in particular, of class $C^{r+1}_*$,
which denotes the Zygmund space
coinciding with the H\"older space $C^{r+1}$
for any positive non-integer $r$
(cf.~proof of Lemma~1.2 in \cite{ht}).
Thus, choosing a non-integer $r$
with $1/2<r\le s$,
we obtain a holomorphic coordinate map
$H$
of class 
$C^{(s+1)-}$
which maps $M$
onto a real-analytic submanifold $\2M$.
The latter fact 
follows by
applying 
the edge-of-the-wedge
theorem as in 
\cite[\S4.3]{KZ}.
Hence $M$
is analytically regularizable
of class 
$C^{(s+1)-}$
at $p$
as desired.

\end{proof}

\section{Sphericity of $3$-dimensional hypesurfaces}

In this section, we 
shall prove \autoref{main-2}
on the
 sphericity 
of real strictly pseudoconvex hypersurfaces $M\subset \C^2$
of any regularity class $C^k$, $k\geq 2$.
We keep all the above notations from Section 2.

In what follows, we will be assuming a regularizing $(0,1)$ vector fields $L$ of class $C^s$, as in \autoref{main-2}, to be {\em basic regularizing}. This is more convenient for our considerations, while it doesn't play a role for our conclusions, as follows from \autoref{cond-indep}. 

\begin{proof}[Proof of necessity] Sphericity of $M$ 
of class $C^s$
in a  neighborhood $U$ of a point $p$ in $M$
 is equivalent to the existence of a CR-diffeomorphism 
 $H$
 of class $C^s$
 of $U$ onto a spherical real-analytic 
 hypersurface $\tilde M$. 
We consider, as above, the image  $M_J$ of $M$ under \eqref{th-section} (which is a $C^{k-1}$ smooth totally real submanifold in $J^{1,n}$). 
Recall that $M_J$ is 
a submanifold in the holomorphically 
degenerate weakly pseudoconvex hypersurface
 $\widehat M=M\times\mathbb{CP}^1\subset J^{1,n}$ 
  of class $C^{k}$.  We next consider the $C^{k-2}$ function $\Phi$ (on either $M$ or $M_J$)
  as in \eqref{Phi-def}
   associating to a point $p\in M$ the $2$-jet of its (formal) Segre variety $Q_p$.

Assume first that $M$ is spherical 
of class $C^s$ at $p\in M$
with $2\le s\le k$.
Then \autoref{maink} (i) implies
 that Condition E holds at $p$ of
 class $C^{s-2}$.
  That is, $\Phi$ must extend  from $M_J$ to the pseudoconvex side of $\widehat M$
  holomorphically and $C^{s-2}$ up to the boundary. 
  Note that, in view of the Lewy extension property (e.g. \cite{ber}), the latter is equivalent to the existence of a $C^{s-2}$  extension of $\Phi$ from $M_J$ to $\widehat M$ as a {\em CR-function} $\varphi$. Furthermore, the analytic image $\tilde M$ of $M$ under the 
  given CR-diffeomorphism must be in addition spherical. 

We now make use of the following remarkable fact, going back to the work of E.~Cartan and A.~Tresse. 
We refer to e.g. \cite{divergence} for more details. 
\begin{proposition}\Label{cubicODE}
A Levi-nondegenerate real-analytic hypersurface $M\subset\CC{2}$ is spherical if and only if the defining function $\Phi$ of its associated ODE  \eqref{wzz} is cubic with respect to the the jet coordinate $\xi:=w'$.  
\end{proposition}

Coming back to our argument, we conclude from
 \autoref{cubicODE} that the associated to $\tilde M$ function 
 $
 \tilde\Phi$, considered already on an open set in $J^{1,n}$, 
 is 
 {\em cubic in the jet variable $\xi$} (see e.g. \cite{divergence} for the discussion of the cubic property and its invariance). In view of the invariance, the extension of the original function $\Phi$ to the pseudoconvex side must be cubic in $\xi$ too. 
 Taking the boundary value $\phi$ of $\Phi$,
 we obtain
\beq\Label{phi}
\varphi
=a_0
+a_1
\xi+a_2\xi^2+a_3\xi^3
\eeq
for some CR-functions $a_0,...,a_3$ on $M$,
which we claim are
of class $C^{s-2}$ --- the same as $\Phi$.
Indeed, the coefficients
of a polynomial can be expressed
in terms of
values at finitely many points
by means of the Lagrange interpolation formula,
which provides the same regularity
as restrictions of $\Phi$ for fixed values of $\xi$.

Extending 
the CR functions $a_0,...,a_3$
holomorphically to the pseudoconvex side
and using the uniqueness of holomorphic extension,
we finally obtain:
\beq\Label{condit}
\Phi
=a_0+a_1\theta+a_2\theta^2+a_3\theta^3,
\eeq
where 
$\theta$ is as in \eqref{th-def},\eqref{th-coor}. 
%
At this point, we  pick a nonzero $(0,1)$ vector field $L$
that is
 {\em regularizing of class $C^s$} 
 (see Definition~\ref{regu}).
 Then $L\th$ is a CR function of class $C^s$
 and
 strict pseudoconvexity implies $L\th\ne 0$.
Hence we can
define
 \beq
 \Label{l-def}
 \mathcal L:=\frac{1}{L\theta}L
 \quad
 \Longrightarrow
 \quad
 \mathcal L\theta=1.
 \eeq

Now, applying $\mathcal L$ once to \eqref{condit} and using \eqref{l-def}, we obtain
\beq\Label{L}
\mathcal L\Phi=a_1+2a_2\theta+3a_3\theta^2,
\eeq
where we took into account the fact that all the $a_j$ are CR-functions. 
Note that the latter implies the regularity condition 
$$
\mathcal L\Phi  \in C^{s-2}(U),\quad L\Phi=L\theta\cdot\mathcal L\Phi\in C^{s-2}(U).
$$
Repeating this argument three more times, we get respectively
\beq\Label{L2}
\mathcal L^2\Phi=2a_2+6a_3\theta,\eeq
then
\beq\Label{L3}
\mathcal L^3\Phi=6a_3,\eeq
and finally:
\beq\Label{spherical}
\mathcal L^4\Phi=0,\quad L^4\Phi=0
\eeq
(see the discussion in the Introduction).
Furthermore, we have the improved regularity properties
\beq\Label{improved}
\mathcal L\Phi,\mathcal L^2\Phi,\mathcal L^3\Phi,L\Phi,L^2\Phi,L^3\Phi\in C^{s-2}(U).
\eeq 
Thus, the differential equation \eqref{spherical} and the regularity conditions \eqref{improved}  give  necessary conditions for the sphericity of $M$, and this proves the necessity in \autoref{main-2}. 
\end{proof}


We are now in the position to prove the sufficiency of our sphericity criterion in low regularity.

\begin{proof}[Proof of sufficiency]
We argue by reversing the above argument for the necessity of \eqref{spherical},\eqref{improved}. 
Let
 $L$ be a
$(0,1)$ vector field on $M$
that is 
{\em regularizing of class $C^s$}, $s\le k$,
i.e.\ $L\th$ is a CR function of class $C^s$ on $M$
(see Definition~\ref{regu}),
and recall that \eqref{improved} holds.
Define the new $(0,1)$ vector field $\6L$
by \eqref{l-def}, hence 
$\6L\th=1$
and by the argument in Remark~\ref{cond-indep},
it follows that
\beq
\Label{main-cond''}
\6L\Phi, \6L^2\Phi, \6L^3\Phi \in C^s,
\quad
\6L^4 \Phi =0.
\eeq
Now
\eqref{main-cond''}
implies that the function $\6L^3\Phi$ is CR of class $C^s$,
hence
there exists a CR function $a_3$  on $M$
satisfying
$$
\6L^3\Phi = 6a_3 \in C^s.
$$
%
%
The latter identity can be writeen as
$$
\6L(\6L^2\Phi) = 6a_3.
$$
Using $\6L\th=1$
along with the CR property of $a_3$, we compute
$$
\6L(\6L^2\Phi - 6 a_3\th) = 6 a_3 - 6 a_3 =0,
$$
hence there exists a CR function $a_2$ satisfying
$$
\6L^2\Phi - 6 a_3\th = 2a_2.
$$
Since $a_3\th \in C^s$,
\eqref{main-cond''} implies $a_2\in C^s$.
Next, we compute
$$
\6L(\6L\Phi  - 3a_3\th^2- 2a_2\th)
=
\6L^2\Phi - 6 a_3\th - 2a_2
=
0,
$$
hence there exists a $C^s$ CR function $a_1$ on $M$ satisfying
$$
\6L\Phi  - 3a_3\th^2- 2a_2\th
=
a_1,
$$
which is again $C^s$
since all other terms are $C^s$.
Finally,
$$
\6L(\Phi - a_3\th^3- a_2\th^2 - a_1\th)
=
\6L\Phi - 3a_3\th^2 - 2a_2\th - a_1
=0,
$$
hence there exists a $C^s$ CR function $a_0$ on $M$ satisfying
$$
\Phi - a_3\th^3- a_2\th^2 - a_1\th = a_0,
$$
which is again of class $C^s$
by the same argument as before.
The latter identity is equivalent to \eqref{condit}
%
that
precisely means that $\Phi$ is the restriction onto $M_J$ of the function $\varphi$, defined as in \eqref{phi},
which is a cubic polynomial with CR coefficients
of class $C^s$.
Then 
by the Lewy extension,
$\phi$ and hence $\Phi$
extend to the pseudoconvex side of $\3M$
 holomorphically and $C^s$ up to the boundary.
 
Now, 
$M$ satisfies
Condition $E$ 
of class $C^s$, and hence 
$M$ is analytically regularizable
of class $C^{(s+1)-}$
 by \autoref{maink}.
More precisely, there exists 
a CR-diffeomorphism of class $C^{(s+1)-}$
between $M$ and a real-analytic strictly pseudoconvex hypersurface
$\2M$.
 
To conclude finally the sphericity of $M$, we use the cubic representation \eqref{phi}
 and its invariance 
  to conclude that
   the ODE associated with $\tilde M$ 
 is also cubic in the respective jet variable, 
 so that by \autoref{cubicODE}, 
 $\tilde M$ is $C^s$ spherical
 and so is $M$, as required. 
\end{proof}

\begin{proof}[Proof of \autoref{main-2-iff} and \autoref{main-2-suff}] The proof is acccomplished by using identically the same argument as that in the above proof of \autoref{main-2}, but employing \autoref{Maink} (more precisely, the last assertion of it) instead of \autoref{maink}. 
\end{proof}

\begin{proof}[Proof of \autoref{zerocurv}] 

For  the proof, we will make use of \autoref{main-2-iff}. Note that, in the regularity $C^k$ with $k=6$, the conditions $\mathcal L\Phi, \mathcal L^2\Phi, \mathcal L^3\Phi\in C^1$ hold automatically (since $\Phi\in C^4$ by its definition). So, it only remains 
to prove
\beq\Label{L4=0}
\mathcal L^4\Phi=0.
\eeq

Recall that we use $\mathcal L$ which is  basic regularizing. 
To prove \eqref{L4=0}, we shall first prove the following claim:

\smallskip

{\em for a $C^6$ strictly pseudoconvex hypersurface $N\subset\CC{2}$, the expression $\mathcal L^4\Phi$ is relatively invariant under local biholomorphisms, i.e. if two coordinate systems are mapped into each other bi a biholomorphism $H$, then the respective expressions $\mathcal L^4\Phi$ for the respective hypersurfaces are related by a nonvanishing factor (actually depending rationally on the $6$-jets of both the map $H$ and the source defining functiong $\rho$, as in \eqref{Phi1}).} 

\smallskip

To prove the claim, we first consider the case of a {\em real-analytic} hypersurface $N$ subject to a biholomorphism $H$. Then $\Phi$ can be considered as holomorphic in a full neighborhood of the origin in $J^{1,n}$. We then observe the identity
\beq\Label{Lxi}
\mathcal L g = \frac{\partial}{\partial \xi}g 
\eeq
which holds at all point on $M_J$ for any function $g$ holomorphic in a neighborhood of the origin in $J^{1,n}$ (actually, one sided neighborhood is sufficient here). The identity follows by subjecting the holomorphic differential $dg$ onto the submanifold $M_J$. Next, \eqref{Lxi} implies that $\mathcal L^4\Phi$ (in the analytic case) coincides with $\frac{\partial^4}{\partial \xi^4}\Phi$. The latter expression is the well known {\em Cartan-Tresse relative invariant} of ODEs (\cite{cartan,tresse}), hence the desired invariance of  $\mathcal L^4\Phi$ follows. Now, in the general $C^6$ smooth case, we observe that the transformation rule for $\mathcal L^4\Phi$ (as a rational function of the $6$-jet of $\rho$)  under a biholomorphism $H$ is the same as in the analytic case, and we finally prove the claim. 
The claim implies that {\em both the identical and the pointwice vanishing of $\mathcal L^4\Phi$ are biholomorphically invariant}. 

On the other end, the CR-curvature is known to be a rational expression in the $6$-jet of the defining function $\rho$ (\cite{cartan,chern}). In the analytic case, when in particular the full Cartan-Tanaka-Chern-Moser theory applies, this expression is known to be (relatively) invariant. Hence it is relatively invaria \TeX also in the $C^6$ case (when it can be formally computed, as a continuous function, but the equivalence theory doesn't apply).

To complete the proof now, we apply the normal form theory in \cite{chern} and conclude that there is a local biholomorphic (in fact bi-polynomial) transformation of $\CC{2}$, mapping a point $q\in M$ to the origin and $M$ into a "normal form to order $6$":
\beq\Label{ord6}
\im w=|z|^2+(c_6z^4\bar z^2+\bar c_6z^2\bar z^4)+o(|z|^6),\quad c_6\in\CC{}.
\eeq
Furthermore, according to \cite{chern}, the coefficient $c_6$ is a nonzero multiple of the CR-curvature at $q$ (this holds in the analytic case, hence by uniqueness in the $C^6$ case too). We get that $c_6=0$. Finally, a direct application of \eqref{Phi1} shows that $\mathcal L^4\Phi$, computed for a hypersurface \eqref{ord6}, vanishes at $0$, hence the respective expression vanishes at $q$. Since $q$ is arbitrary, this prove \eqref{L4=0} and the theorem.

\end{proof}

\end{document}